\documentclass[aap,seceqn,citesort,MSNbibl,dvips]{arximspdf}

\doi{10.1214/09-AAP639}
\volume{20}
\issue{3}
\pubyear{2010}
\firstpage{785}
\lastpage{805}

\makeatletter

\newtheorem{prop}{Proposition}[section]
\newtheorem{theorem}[prop]{Theorem}

\newproclaim{rem}[prop]{Remark}
\newproclaim{examp}[prop]{Example}

\newtheorem{lem}[prop]{Lemma}

\makeatother

\begin{document}
\begin{frontmatter}

\title{Intermediate range migration in the two-dimensional
stepping stone model}
\runtitle{Stepping stone model}

\begin{aug}
\author[A]{\fnms{J. Theodore} \snm{Cox}\corref{}\thanksref{T1}\ead[label=e1]{jtcox@syr.edu}}
\runauthor{J. T. Cox}
\affiliation{Syracuse University}
\address[A]{Mathematics Department\\
Syracuse University\\
Syracuse, New York 13244\\
USA\\
\printead{e1}} 
\end{aug}

\thankstext{T1}{Supported in part by NSF Grant 0803517.}

\received{\smonth{7} \syear{2009}}
\revised{\smonth{1} \syear{2010}}

%
\begin{abstract}
We consider the stepping stone model on the torus of side $L$ in
$\mathbb{Z}^2$
in the limit $L\to\infty$, and study the time it takes two lineages
tracing backward in time to coalesce. Our work fills a gap between the
finite range migration case of [\textit{Ann. Appl. Probab.} \textbf{15}
(2005) 671--699] and the long range case of [\textit{Genetics}
\textbf{172} (2006) 701--708], where the migration range is a positive
fraction of $L$. We obtain limit theorems for the intermediate case,
and verify a conjecture in [\textit{Probability Models for DNA Sequence
Evolution} (2008) Springer] that the model is homogeneously mixing if
and only if the migration range is of larger order than $(\log
L)^{1/2}$.
\end{abstract}

%
\begin{keyword}[class=AMS]
\kwd[Primary ]{60K35}
\kwd{60G50}
\kwd{92D10}
\kwd[; secondary ]{82C41}.
\end{keyword}
\begin{keyword}
\kwd{Stepping stone model}
\kwd{torus random walk}
\kwd{hitting times}
\kwd{coalescence times}.
\end{keyword}

\end{frontmatter}

\section{Introduction}

The subject of this paper is the stepping stone model of
population genetics, and in particular the contrast between
recent results of \cite{MW06} and \cite{ZCD} in the
two-dimensional setting. There is a vast literature on the
many variants of the stepping stone model dating back to the
seminal work of Mal\`{e}cot \cite{Ma49} and Kimura
\cite{Ki53}. (A few sources for background and references
are \cite{CG87,Du08,Na89} and \cite{Wi04}.) We will
begin by describing the version of the model we consider
here, generally following the setup in \cite{ZCD}.

Let $\mathbb{Z}^2$ be the two-dimensional integer lattice, and fix
$\nu>0$ and $q\dvtx\mathbb{Z}^2\to[0,1]$ with $q(0)=0$ and $\sum_x
q(x)=1$. We suppose that at each site $x$ in
\[
\mathbb{T}_L=(-L/2,L/2]^2\cap\mathbb{Z}^2
\]
there is a colony of $2N$ haploid individuals. We think of
$\mathbb{T}_L$ as a torus, and assume a continuous-time Moran model
of reproduction. In this model, a given individual in colony
$x$ dies at rate one, independently of all other
individuals, and is replaced by a copy of an individual chosen
at random from the same colony with probability $1-\nu$ or
colony $y$ with probability $\nu q(y-x)$ computed modulo
$L$. In this way, we treat $\mathbb{T}_L$ as a torus. The
genealogical structure of a sample of $n$ individuals is
determined by tracing their lineages backward in time.

We will focus on the case of $n=2$ lineages, where one is
interested in $T_0$, the time it takes the lineages
to enter the same colony, and $t_0$, the time to coalescence
of the lineages. There are many limit
theorems for $T_0$ and $t_0$ in the literature. (A~small
sampling can be found \cite{FOW,C89,WH98,CD02,Wi04,MW06}
and \cite{ZCD}.) One
may allow $N\to\infty$, $\nu\to0$ and $q$ to vary as
$L\to\infty$. To understand the asymptotic behavior of $t_0$,
one must first understand the behavior of $T_0$ so
we will concentrate on the latter. Furthermore, the question
we want to consider is already of interest in the simplest
case of one individual per colony, so we will assume from
now on that $\nu=2N=1$, but allow $q$ to vary.

The \textit{meanfield} or \textit{homogeneous mixing} case is
obtained by taking $q$ to be uniform over
$\mathbb{T}_L\setminus\{0\}$. Suppose the two lineages start at
$0,x\in\mathbb{T}_L$, $x\ne0$. The law of $T_0$
is exponential with mean $(L^2-1)/2$ and is independent of
$x$, and so $T_0/L^2$ converges in law, uniformly in
$x\ne0$, to the exponential distribution with mean~$1/2$.
Matsen and Wakeley show in \cite{MW06} that the same
limiting behavior of $T_0/L^2$ holds uniformly in $x\ne0$
assuming that $q$ is uniform on only a positive fraction
of the torus. By contrast, if $q$ is kept fixed as
$L\to\infty$, then the right normalization for $T_0$ is
$L^2\log L$, and the limiting law depends on the
starting positions $0,x$. (See \cite{C89,CD02} and
\cite{ZCD} for results of this type.) The purpose of this paper
is to fill the gap between these two situations.

Following two lineages backward in time amounts to following
two random walks until they meet. The
difference between the lineage locations is also a random
walk, and $T_0$ is just the time it takes this difference walk
to hit 0. On account of this, we will now focus on the following
random walk setting.
For $k>0$, let
\[
\Lambda_k= [-k/2,k/2]^2\cap\mathbb{Z}^2
\]
and for any $A\subset\mathbb{R}^2$ let
\[
A'=A\setminus\{0\}.
\]
For $r>0$, let $B(r) = \{x\in\mathbb{R}^2\dvtx\|x\|_\infty\le
r\}$. Let $M_1,M_2,\ldots$ be a sequence of positive integers
and assume that $q_{M_L}\dvtx\mathbb{Z}^2\to[0,1]$ satisfies\vspace*{22pt}
\begin{enumerate}[(P0)]
\item[(P0)]\hypertarget{equationP0}
\mbox{}
\begin{eqnarray*}
\\[-60pt]
q_{M_L}(x)&=&0 \qquad\mbox{for }x\notin\Lambda'_{M_L},\nonumber\\
\sum_x q_{M_L}(x)&=&1,
q_{M_L} \mbox{ is symmetric}\quad\mbox{and}\\
\sigma_{M_L}^2 &\equiv& \sum_x x^2_1 q_{M_L}(x)
= \sum_x x^2_2q_{M_L}(x) >0.\nonumber
\end{eqnarray*}
\end{enumerate}
The uniform distributions $u_{M_L}(x)
= 1_{\Lambda'_{M_L}}(x)/|\Lambda'_{M_L}|$ clearly satisfy
\hyperlink{equationP0}{(P0)}.

Let\vspace*{1pt} $Y^L_t$ be a rate one random walk on $\mathbb{Z}^2$ with jump
distribution $q_{M_L}$, and let~$X^L_t$ be the corresponding walk
on $\mathbb{T}_L$ viewed as a torus. Given $Y^L_t$ we
construct $X^L_t$ by setting $X^L_t=Y_t\operatorname{mod} L$.
Let $H_{L}$ be the hitting time for $X^L_t$ of the origin,
\[
H_{L} = \inf\{t\ge0\dvtx X^L_t=0\}.
\]
Then $H_L$ has the same law as $2T_0$, so we will study
$H_L$. Let $P_x$ and $E_x$ denote probability law and
expectation for the walk starting at $x$.

With the above notation, the Matsen and Wakely result is as
follows. Fix $0<c<1$ and let $M_L=cL$ and
$q_{M_L}=u_{M_L}$. Then
as $L\to\infty$,
%
%
\begin{equation}\label{mfexplaw}
H_L/L^2 \Rightarrow\mathcal{E}(1) \qquad\mbox{uniformly in }
X^L_0\in\mathbb{T}_L',
\end{equation}
where $\Rightarrow$ indicates the law of the left-hand side converges
weakly to the distribution on right-hand side, and $\mathcal{E}(\beta)$
is the exponential distribution with mean $\beta$. On the
other hand, if $M_L\equiv M$ is fixed, so there is a single
jump distribution $q_M$, then by Theorem 1 of \cite{ZCD},
if $0<\alpha< 1$ and $|X^L_0|\approx L^{\alpha}$ as
$L\to\infty$, then
%
%
\begin{equation}\label{myexplaw}
\frac{H_L}{L^2\log L} \Rightarrow(1-\alpha)\delta_0 + \alpha
\mathcal{E}(1/\pi\sigma^2_M).
\end{equation}
Here, $x_L\approx L^\alpha$ means $x_L\in\mathbb{T}_{L^\alpha\log
L}\setminus\mathbb{T}_{L^\alpha/\log L}$.

It seems clear that the homogeneous mixing behavior of
(\ref{mfexplaw}) should hold if $M_L\to\infty$ at a
sufficiently fast rate, and it is natural to ask what this
rate might be. Durrett (see Section 5.6 and Theorem 5.18 of
\cite{Du08}) conjectured that it should be quite slow, only of
greater order than $\sqrt{\log L}$ as $L\to\infty$, meaning
that (\ref{mfexplaw}) should hold exactly when $M_L/\sqrt{\log
L}\to\infty$. We verify this conjecture for a large class of
jump distributions in Theorems \ref{thm:rick} and
\ref{thm:main} below, and obtain a slightly improved version
of (\ref{myexplaw}) when $M_L=O(\sqrt{\log L})$. The proof of
(\ref{mfexplaw}) in \cite{MW06} makes use of Markov chain
techniques from \cite{A89} and \cite{DS91}. The proof of
(\ref{myexplaw}) relies heavily on local central limit theorem
estimates for $P_0(Y^L_t=0)$ to then estimate $P_0(X^L_t=x)$
[for use in (\ref{FL}) below]. Here, we will use a more direct
Fourier-type approach that seems simpler, and works for both
(\ref{mfexplaw}) and (\ref{myexplaw}) as well.

For a jump distribution $q_M$, define the characteristic function
\[
\phi_M(\theta) = \sum_{x\in\Lambda_{M}}e^{i\theta x}q_M(x)
,\qquad
\theta\in\mathbb{R}^2,
\]
where $\theta x=\theta\cdot x$. We will
assume that the jump distributions $q_{M_L}$ have
characteristic functions $\phi_{M_L}$ which satisfy the
conditions \hyperlink{equationP1}{(P1)}--\hyperlink{equationP3}{(P3)} listed below. These conditions are satisfied
for the uniform distributions $u_{M_L}$ (see the Appendix of
\cite{CDP00}, where $M^2$ in \hyperlink{equationP2}{(P2)} there should be $M$).
Proposition \ref{prop:suff} below shows that they are satisfied
in some generality. Note that the symmetry condition in \hyperlink{equationP0}{(P0)} implies
each $\phi_M$ is real-valued. The conditions we need are the following.
\begin{enumerate}[(P1)]
\item[(P1)]\hypertarget{equationP1}
There is a $\sigma^2>0$ such
that for all $\varepsilon>0$ there exists $\delta>0$ such that for all
large $L$,
\[
\dfrac{1-\phi_{M_L}(\theta)}{\sigma^2M_L^2|\theta|^2/2} \in
(1-\varepsilon,1+\varepsilon) \qquad\mbox{for all }
\theta\in B'(\delta/M_L).
\]
\item[(P2)]\hypertarget{equationP2}
For all $\delta>0$, there exists
$\delta'>0$ and $\zeta>0$ such that for all large
$L$,
\[
1-\phi_{M_L}(\theta) > \zeta\qquad\mbox{for
all }\theta\in B(\delta')
\setminus B(\delta/M_L).
\]
\item[(P3)]\hypertarget{equationP3}
For all $\varepsilon>0$ and $a>0$,
\[
|\phi_{M_L}(\theta)| < \varepsilon\qquad\mbox{for all }\theta\in B(\pi)
\setminus B(a)
\mbox{ and all large }L.
\]
\end{enumerate}
\begin{prop}\label{prop:suff} Let $f$ be a positive,
continuous function on $B(1/2)$ such that
$f(x_1,x_2)=f(x_2,x_1) = f(-x_1,x_2)$. Define $c_M>0$ and
$q_M(x) = c_Mf(x/M)u_M(x)$ so that $\sum_xq_M(x)=1$. Then
for any $M_L\to\infty$ as $L\to\infty$, the corresponding
sequence of characteristic functions $\phi_{M_L}$
satisfies properties \textup{\hyperlink{equationP1}{(P1)}}--\textup{\hyperlink{equationP3}{(P3)}}.
\end{prop}

In addition to \hyperlink{equationP1}{(P1)}--\hyperlink{equationP3}{(P3)}, we impose the mild regularity
condition\vspace*{9pt}
\begin{enumerate}[(P4)]
\item[(P4)]\hypertarget{equationP4}
\mbox{}
\begin{eqnarray*}
\\[-43pt]
\lim_{L\to\infty}\frac{M_L^2}{\log L} =
\rho\in[0,\infty].
\end{eqnarray*}
\end{enumerate}

Our first result shows that homogeneous
mixing occurs if
$M^2_L/\log L\to\infty$.
\begin{theorem}\label{thm:rick}
Assume conditions
\textup{\hyperlink{equationP0}{(P0)}}--\textup{\hyperlink{equationP4}{(P4)}}
hold with $\rho=\infty$. Then for all $\lambda>0$,
%
%
\begin{equation}\label{eqn:exp-law-1}
{\lim_{L\to\infty}\sup_{x\in\mathbb{T}_L'}} |
E_x(e^{-\lambda H_L/L^2}) - ( 1 + \lambda)^{-1}
| = 0
\end{equation}
and
%
%
\begin{equation}\label{eqn:mean1}
{\lim_{L\to\infty}\sup_{x\in\mathbb{T}_L'} }|
E_x(H_L/L^2) - 1 | = 0.
\end{equation}
\end{theorem}

Our next result shows that homogeneous mixing does \textit{not} occur if
$\rho<\infty$, and that $H_L$ can grow at any rate
between $L^2$ and $L^2\log L$. We will use the following notation.
For $v>0$, define
\[
\mathcal{A}_L(\alpha,v) =
\cases{
\mathbb{T}'_v, &\quad if $\alpha=0$,\cr
\mathbb{T}_{L^\alpha v} \setminus
\mathbb{T}_{L^\alpha/v}, &\quad if $0<\alpha<1$,\cr
\mathbb{T}_L\setminus\mathbb{T}_{L/v},&\quad if $\alpha=1$,}
\]
and let
%
%
\begin{equation}\label{betadef}
t_L=\frac{\log L}{M^2_L} \quad\mbox{and}\quad
\beta= \rho+\frac{1}{\pi\sigma^2}.
\end{equation}
\begin{theorem}\label{thm:main}
Assume $M_L\to\infty$ and the conditions
\textup{\hyperlink{equationP0}{(P0)}}--\textup{\hyperlink{equationP4}{(P4)}}
hold with $\rho<\infty$. Fix $0\le\alpha\le1$ and $k>0$, and put
$v_L=(\log L)^k$. Then for all $\lambda>0$,
%
%
\begin{equation}\label{eqn:exp-law}\quad
{\lim_{L\to\infty}\sup_{x\in\mathcal{A}_L(\alpha,v_L)}} |
E_x(e^{-\lambda H_L/L^2t_L}) - [(1-\alpha') + \alpha'( 1 +
\beta\lambda)^{-1} ]
| = 0
\end{equation}
where $\alpha' = (\alpha+\rho\pi\sigma^2)/(1+\rho\pi\sigma^2)$.
Furthermore,
%
%
\begin{equation}\label{eqn:mean}
{\lim_{L\to\infty}\sup_{x\in\mathcal{A}_L(\alpha,v_L)}} |
E_x(H_L/L^2t_L) - \alpha'\beta| =0.
\end{equation}
\end{theorem}
\begin{rem}If we set $\rho=0$ in (\ref{eqn:exp-law}), then we
recover the form (\ref{myexplaw}). The proof of
(\ref{eqn:exp-law}) is easily adapted to handle the case of
a fixed $q_M$ satisfying \hyperlink{equationP0}{(P0)}, providing a slight
strengthening of (\ref{myexplaw}). One can also see that
(\ref{eqn:exp-law}) is consistent with
(\ref{eqn:exp-law-1}) by setting $t_L\equiv
1$, rephrasing (\ref{eqn:exp-law}) appropriately, and then
setting \mbox{$\rho=\infty$}.
\end{rem}

A one-dimensional stepping stone model was considered in
\cite{DR08}, where exponential limit laws for $H_L$ were
obtained under rather general assumptions on the
jump distributions. We will not state their results, but note that in
analogy with Theorem~3 there, one might hope in our
two-dimensional setting that with $M^2_L=\log L$ some
version of (\ref{myexplaw}) would hold with
\hyperlink{equationP1}{(P1)}--\hyperlink{equationP3}{(P3)} replaced by the simpler conditions
%
%
\begin{eqnarray}\label{insuffconds}
\mbox{(i) }&&\lim_{L\to\infty}\sigma^2_{M_L}/M^2_L=\sigma^2\quad
\mbox{and } \nonumber\\[-8pt]\\[-8pt]
\mbox{(ii) }&&\mbox{for some }c>0,\qquad q_{M_L} \ge c u_{M_L}.\nonumber
\end{eqnarray}
More precisely, the desired result would be that
(\ref{insuffconds}) implies $H_L/L^2\Rightarrow\mathcal{E}(\beta
_0)$ for
$X^L_0$ large, where the limiting mean $\beta_0$ depends only
on $\sigma^2$ and $c$. This is not the case, as the following
example shows.
\begin{examp}\label{myexample}
Fix $0<c<1$ and $q_0\dvtx\mathbb{Z}^2\to[0,1]$ satisfying
\hyperlink{equationP0}{(P0)} for some fixed $M_0$, and let $\hat
q_0(\theta) = \sum_x q_0(x)e^{i\theta x}$. Put $q_{M_L}(x)
= c u_{M_L}(x) + (1-c)q_{0}(x)$, assume that
$\lim_{L\to\infty}M^2_L/\log L= 1$, and define
%
%
\begin{equation}\label{beta0}
\beta_0 = \frac{12}{c\pi} +
\frac1{(2\pi)^2}\int_{B(\pi)}\frac{d\theta}{1-(1-c)\hat
q_0(\theta)}.
\end{equation}
Then $q_{M_L}$ satisfies (\ref{insuffconds}) with
$\sigma^2=c/12$. If $L/\log L < \ell_L<L$,
then for all $\lambda>0$,
%
%
\begin{equation}\label{eqn:mix}
{\sup_{x\in\mathbb{T}_L\setminus\mathbb{T}_{\ell_L}}} |
E_x(e^{-\lambda H_L/L^2}) - ( 1 + \beta_0\lambda)^{-1}
| \to0\qquad
\mbox{as }L\to\infty.
\end{equation}
\end{examp}
\begin{rem} The influence of the short range jumps is
reflected in the dependence of $\beta_0$ on $\hat
q_0$. Other mixtures of jump distributions could also be
considered, e.g., $\sum_i c_i u_{M^i_L}$ where
$M^1_L,M^2_L,\ldots$ tend to infinity at different rates.
\end{rem}

The proofs in \cite{C89,CD02} and \cite{ZCD} for the
fixed jump distribution case use the fact that $X^L_t$
becomes uniformly distributed over the torus by times of
larger order than~$L^2$. The analogous fact in our setting
is given below, it will be used in the proof
of~(\ref{eqn:mean}).
\begin{theorem}\label{thm:uniform} Assume
\textup{\hyperlink{equationP1}{(P1)}}--\textup{\hyperlink{equationP4}{(P4)}}
hold. If $s_L/[ (L^2/M_L^2)\vee\log
L]\to\infty$ as $L\to\infty$, then
%
%
\begin{equation}\label{eqn:uniform}
\lim_{L\to\infty}\sup_{t\ge s_L}\sup_{x\in\mathbb{T}_L} L^2|
P_0(X^L_t=x)-L^{-2}| = 0.
\end{equation}
\end{theorem}

Returning to the stepping stone model, we could now
consider the genealogy of a sample of $n>2$ individuals.
Let $\zeta^L_t$ be a system of rate one coalescing random
walks on $\mathbb{T}_L$ with jump distribution
$q_{M_L}$. If we consider
lineages starting at $x_i\in\mathbb{T}_L,1\le i\le n$, and put
$\zeta^L_0=\{x_1,\ldots, x_n \}$, then $|\zeta^L_t|$ is the
number of distinct lineages left at time $t$. Under the assumptions of
Theorem \ref{thm:main}, and assuming $|x_i-x_j|\ge
L/\log L$ for $i\ne j$, the analog of Theorem 2 of
\cite{ZCD} would be
%
%
\begin{equation}\label{nCRW}
\lim_{L\to\infty}P(|\zeta^L_{sL^2t_L}| = k)= P(D_{\pi\sigma
^2s}=k),\qquad k=1,\ldots, n,
\end{equation}
where $D_t$ is the pure death
process on the positive integers which makes transition $k\to
k-1$ at rate ${k\choose2}$. In fact, the genealogy of the
lineages (on this time scale) converges to the genealogy described
by Kingman's
coalescent (see \cite{Kn82}). We will not pursue these matters here,
since with the results developed the methods of \cite{C89,CD02}
and \cite{ZCD} could be adapted to prove such
limit laws.

The outline of the rest of the paper is as follows. In Section
\ref{sec2}, we develop some simple Fourier analytic tools.
Proposition \ref{prop:suff} is proved
in Section \ref{sec3}, Theorem \ref{thm:uniform} is proved in Section
\ref{sec4}, Theorem \ref{thm:rick} is proved in Section \ref{sec5}, and
Theorem \ref{thm:main} is proved in Section \ref{sec6}. Finally,
we verify the claims for Example \ref{myexample} in
Section \ref{sec7}. For simplicity, we will assume throughout the rest of
the paper that $L,M,M_L,\ldots$ are positive even integers.

\section{Preliminaries}\label{sec2} For a jump distribution $q_M$
satisfying \hyperlink{equationP0}{(P0)} with characteristic
function $\phi_M$, define the transforms
%
%
\begin{eqnarray}\label{phit}
\phi_M^{t}(\theta) &=& E_0(e^{i\theta Y^L_t}) =
\exp\bigl(-t\bigl(1-\phi_M(\theta)\bigr)\bigr),\nonumber\\
F_L(x,\lambda)&=&E_x(e^{-\lambda H_L})\quad \mbox{and }\\
G_L(x,\lambda)&=&\int_0^\infty e^{-\lambda s}
P_x(X^L_s=0)\,ds,\nonumber
\end{eqnarray}
where $\theta\in\mathbb{R}^2$, $t\ge0, x\in\mathbb{T}_L$ and
$\lambda\ge0$.
The reason for our interest in $G_L(x,\lambda)$ is the formula
%
%
\begin{equation}\label{FL}
F_L(x,\lambda) = \frac{G_L(x,\lambda)}{G_L(0,\lambda)},
\end{equation}
a simple consequence of the strong Markov property.
We will also make use of the well-known Fourier inversion
formula
%
%
\begin{equation}\label{ptL}
P_0(X^L_t=x) = \frac1 {L^2} \sum_{y\in\mathbb{T}_L} \phi
_M^{t}(2\pi y/L)
e^{2\pi i x y/L},\qquad x\in\mathbb{T}_L,
\end{equation}
from which it is easy to derive
%
%
\begin{equation}\label{GL}
G_L(x,\lambda) = \frac1 {L^2} \sum_{y\in\mathbb{T}_L}
\frac{e^{2\pi i x y/L}}{1+\lambda- \phi_M(2\pi i y/L)}.
\end{equation}

In order to obtain useful bounds on the above, we
will need to estimate sums of complex exponentials over various
regions, including
\[
D_k = \{x\in\mathbb{Z}^2\dvtx|x|\le k/2\},
\]
where $|x|=\|x\|_2$.
\begin{lem}\label{lem:phiest} \textup{(a)} For $K\ge1$ and $\theta\in B(\pi)$,
%
%
\begin{eqnarray}\label{phiest}
\biggl|\sum_{x\in\mathbb{T}_K}
e^{i\theta x} \biggr| &\le&
4(K+1)(1+\|\theta\|^{-1}_\infty)\quad\mbox{and}\nonumber\\[-8pt]\\[-8pt]
\biggl|\sum_{x\in D_K}
e^{i\theta x} \biggr| &\le& 4(K+1)\|\theta\|^{-1}_{\infty}.\nonumber
\end{eqnarray}

\textup{(b)} There is a constant $C_0$ such that for all $J\ge1$ and
$\theta\in B'(\pi)$,
%
%
\begin{equation}\label{DKJbound}
\sup_{K>J} \biggl| \sum_{y\in D_K\setminus D_J}\frac{e^{i\theta y}}{|y|^2}
\biggr| \le\frac{C_0}{1\wedge(J\|\theta\|_\infty)}.
\end{equation}

\textup{(c)}
%
%
\begin{equation}\label{twopi}
\lim_{K\to\infty}\frac{1}{\log K} \sum_{y\in\mathbb{T}'_{K}}
\frac{1}{|y|^2} =
\lim_{K\to\infty}\frac{1}{\log K} \sum_{y\in D'_{K}}
\frac{1}{|y|^2}
= 2\pi.
\end{equation}
\end{lem}
\begin{pf}
Combining the\vspace*{1pt} two elementary facts
$\sin u\ge u/2$ for
$|u|\le\pi/2$ and
$\sum_{j=-k}^k e^{iju} = {\sin((k+\frac12)u)}/{\sin
\frac u2}$ for any positive integer $k$ and real $u$ we obtain
\[
\Biggl|\sum_{j=-k}^k e^{iju}\Biggr| \le4/|u| \qquad\mbox{for all }
k\in\mathbb{Z}^+, u\in B(\pi).
\]
Consequently,
\[
\biggl|\sum_{x\in\Lambda_K} e^{i\theta x} \biggr|\le
\sum_{k=- K/2}^{ K/2}
\Biggl|\sum_{j=-K/2}^{ K/2} e^{i\theta_2 j}
\Biggr| \le\frac{4(K+1)}{|\theta_2|}.
\]
This bound holds with $\theta_1$ replacing $\theta_2$, and
therefore
%
%
\begin{equation}\label{LamLest}
\biggl|\sum_{x\in\Lambda_K}e^{i\theta x} \biggr| \le
4(K+1)\|\theta\|^{-1}_\infty
\qquad\mbox{for all }\theta\in B(\pi).
\end{equation}
The first bound in (\ref{phiest}) follows from this inequality
and the fact that $|\Lambda_K\setminus\mathbb{T}_K|= 2K+1$.
The second bound in (\ref{phiest}) is derived using the
argument for (\ref{LamLest}).

For (b), if $1\le k\le|y|\le
k+1$, then
\[
0 \le\frac1{k^2} - \frac1{|y|^2} \le
\frac1{k^2} - \frac1{(k+1)^2}
\le\frac{6}{|y|^3}.
\]
Let $\gamma_k(\theta) = k^{-2}\sum_{x\in
D_k}e^{i\theta x}$ and
$C=6\sum_{y\in\mathbb{Z}^2\setminus\{0\}}|y|^{-3}<\infty$. Then
\[
\biggl|\sum_{y\in D_K\setminus D_J}\frac{e^{i\theta y}}{|y|^2}
\biggr| \le C +
\Biggl|\sum_{k=J}^{K-1}
\sum_{y\in D_{k+1}\setminus D_k}\frac{e^{i\theta
y}}{k^2}\Biggr|.
\]
We can rewrite the sum on the right-hand side above,
obtaining
\begin{eqnarray*}
\Biggl|\sum_{k=J}^{K-1}
\sum_{y\in D_{k+1}\setminus D_k}\frac{e^{i\theta
y}}{k^2} \Biggr|
&=& \Biggl|
\sum_{k=J}^{K-1}
\biggl(\frac{(k+1)^2}{k^2}\gamma_{k+1}(\theta)-\gamma_k(\theta) \biggr)
\Biggr|\\
&=&
\Biggl|\sum_{k=J}^{K-1}
\biggl(\gamma_{k+1}(\theta)-\gamma_k(\theta) +
\frac{2k+1}{k^2} \gamma_{k+1}(\theta) \biggr) \Biggr|\\
&\le&
|\gamma_{K}(\theta)|+|\gamma_J(\theta)| +
3\sum_{k=J}^{K-1}\frac{|\gamma_{k+1}(\theta)|}{k}.
\end{eqnarray*}
By the bound (\ref{phiest}),
\[
3\sum_{k=J}^{K-1}\frac{|\gamma_{k+1}(\theta)|}{k} \le
\frac{18}{\|\theta\|_\infty}
\sum_{k=J}^{K-1}\frac{1}{k(k+1)}
\le\frac{18}{J\|\theta\|_\infty}.
\]
Making use of the trivial bound $|\gamma_{k}(\theta)|\le(k+1)^2/k^2
\le
4$ for $|\gamma_K(\theta)|$ and $|\gamma_J(\theta)|$, we
therefore have
\[
\biggl| \sum_{y\in D_{K}\setminus
D_J}\frac{e^{i\theta}}{|y|^2} \biggr| \le
C + 8 + \frac{18}{J\|\theta\|},
\]
proving (\ref{DKJbound}).

The second limit in (c) follows from a simple comparison with an
integral. The first follows from a second comparison
showing that
%
%
\begin{equation}\label{2pilog2}
\lim_{K\to\infty}\sum_{D_{2K}\setminus D_K}\frac{1}{|y|^2} = 2\pi
\log2.
\end{equation}
\upqed\end{pf}

We close this section by recording the fact
%
%
\begin{equation}\label{zero}
\sum_{y\in\mathbb{T}_L} e^{2\pi i x y/L} = 0 \qquad\mbox{for all }x\in
\mathbb{T}'_L.
\end{equation}

\section[Proof of Proposition 1.1]{Proof of Proposition \protect\ref{prop:suff}}\label{sec3}
Throughout this section, we will write $M$ for~$M_L$.
It is straightforward to check that the assumptions of
Proposition \ref{prop:suff} imply the following.
As $M\to\infty$:
\begin{longlist}
\item$c_M\to c_0=1/\int_{B(1/2)}f(x) \,dx$,
\item$\sigma^2_M/M^2\to
\sigma^2=c_0\int_{B(1/2)}x_1^2 f(x)\, dx$ and
\item$\phi_M(\theta/M) \to\tilde\phi(\theta) =
c_0\int_{B(1/2)}e^{i\theta x}f(x) \,dx$, $\theta\in B(\pi)$.
\end{longlist}

Let $Z_M$ have distribution
$q_M$. By a standard inequality (see (2.3.6) in \cite{Du05}) and
the fact that $|Z_M|\le M/2$,
%
%
\begin{eqnarray}\label{cfbound}
\bigl|1-\phi_M(\theta)-(\sigma_M^2|\theta|^2/2)\bigr|
&\le& E\bigl( (|\theta Z_M|^3/6)\wedge|\theta
Z_M|^2\bigr)\nonumber\\[-8pt]\\[-8pt]
&\le& \frac{|\theta|^2 M^2}{4} \bigl((|\theta|M/12)
\wedge1\bigr).\nonumber
\end{eqnarray}
Using (ii), this implies that for any $\delta>0$,
\begin{eqnarray*}
\sup_{\theta\in B'(\delta/M)} \biggl|\frac{1-\phi_M(\theta)}{
\sigma_M^2|\theta|^2/2} -1 \biggr| &\le&
\frac{2 M^2}{4\sigma_M^2} \bigl( (\delta/12)\wedge1\bigr)\\
&\to&\frac{1}{2\sigma^2} \bigl( (\delta/12)\wedge1\bigr)\qquad
\mbox{as } M\to\infty.
\end{eqnarray*}
Using (ii) again, this is enough to establish \hyperlink{equationP1}{(P1)}.

Fix $\varepsilon>0$ and put $\bar c=\sup\{c_M\}$. We will prove that
there exists a finite constant $A$ depending on $\varepsilon$ such that
%
%
\begin{equation}\label{p3}
\limsup_{M\to\infty} \sup_{\theta\in B(\pi)\setminus B(A/M)}
\sum_{x\in\Lambda_M}e^{i\theta x}q_M(x)
\le\varepsilon\bar c (1+20\|f\|_\infty),
\end{equation}
which is stronger than \hyperlink{equationP3}{(P3)}. First, we replace the
sum over $\Lambda_M$ with one over $\mathbb{T}_M$ at the cost of a
small error,
%
%
\begin{equation}\label{c1}
\biggl|
\sum_{x\in\Lambda_M}e^{i\theta x}q_M(x) -
\frac{c_M}{|\Lambda'_M|}\sum_{x\in\mathbb{T}_M}e^{i\theta x}f(x/M)
\biggr| \le\frac{(2M+2)\bar c\|f\|_\infty}{|\Lambda'_M|}.
\end{equation}
The idea now is to break the sum over $\mathbb{T}_M$ into sums over
disjoint translates of $\mathbb{T}_K$, where $K<M$ is chosen so that
$f(x/M)$ is essentially constant on the translates, and
then apply (\ref{phiest}).

To do this, let $\Gamma_{M,K}=\{z\in K\mathbb{Z}^2\dvtx z+\mathbb
{T}_K\subset
\mathbb{T}_M\}$, and choose $\varepsilon'\in(0,\varepsilon)$ small
enough so
that $|f(x)-f(x')|<\varepsilon$ if $\|x-x'\|_\infty<\varepsilon'$. Choose
$A$ large enough so that
$A\varepsilon' > \varepsilon^{-1}$ and suppose $\|\theta\|_\infty>A/M$.
Since $|\mathbb{T}_M\setminus\bigcup_{z\in\Gamma_{M,K}}(z+\mathbb
{T}_K)| \le
4MK$,
%
%
\begin{equation}\label{c2}\hspace*{28pt}
\biggl|\sum_{x\in\mathbb{T}_M}e^{i\theta x}f(x/M) -
\sum_{z\in\Gamma_{M,K}}\sum_{x\in\mathbb{T}_K}e^{i\theta(z+x)}
f\bigl((z+x)/M\bigr) \biggr|\le4\|f\|_\infty KM.
\end{equation}
For large $M$, we can choose $K$ to satisfy
$\varepsilon'/2<K/M<\varepsilon'$. By our
choice of $\varepsilon'$, $|f((z+x)/M)-f(z/M)|<\varepsilon$ for all
$z\in\Gamma_{M,K}$ and $x\in\mathbb{T}_K$. Applying this bound
gives
%
%
\begin{eqnarray}\label{c3}
&&\biggl|\sum_{z\in\Gamma_{M,K}}\sum_{x\in\mathbb{T}_K}e^{i\theta(z+x)}
f\bigl((z+x)/M\bigr) \nonumber\\[-8pt]\\[-8pt]
&&\hspace*{14.2pt}{} -
\sum_{z\in\Gamma_{M,K}}e^{i\theta z}f(z/M)\sum_{x\in\mathbb{T}
_K}e^{i\theta x}
\biggr|\le\varepsilon M^2.\nonumber
\end{eqnarray}
By (\ref{phiest}) and the bound $|\Gamma_{M,K}|\le
M^2/K^2$,
%
%
\begin{eqnarray}\label{c4}
\biggl|\sum_{z\in\Gamma_{M,K}}e^{i\theta z}f(z/M)
\sum_{x\in\mathbb{T}_K}e^{i\theta x} \biggr|
&\le& \frac{M^2\|f\|_\infty}{K^2}
4(K+1)(1+\|\theta\|^{-1})
\nonumber\\[-8pt]\\[-8pt]
&\le&\frac{8M^2\|f\|_\infty}{K}(1+\|\theta\|^{-1}_\infty)
.\nonumber
\end{eqnarray}
By combining (\ref{c1})--(\ref{c4}), the bounds
$\varepsilon'/2<K/M<\varepsilon'$ and then using $\|\theta\|_\infty>A/M$,
we obtain
\begin{eqnarray*}
|\phi_M(\theta)| &\le&
\frac{\bar c}{|\Lambda'_M|}
\biggl[(2M+2)\|f\|_\infty
+ 4KM\|f\|_\infty\\
&&\hspace*{30pt}{} +
\varepsilon|M^2| +
\frac{8M^2}{K}\|f\|_\infty(1+\|\theta\|^{-1}_\infty)
\biggr]\\
&\le&\frac{\bar c}{|\Lambda'_M|}
\bigl[
\varepsilon|M^2| + \|f\|_\infty
\bigl((2M+2) + 4\varepsilon' M^2 +
16(M/\varepsilon')(1+M/A)
\bigr) \bigr]\\
&\to&\bar c \biggl[\varepsilon+ 4\varepsilon' \|f\|_\infty+
\frac{16\|f\|_\infty}{\varepsilon' A} \biggr]
\qquad\mbox{as }M\to\infty.
\end{eqnarray*}
Since $A\varepsilon'>\varepsilon^{-1}$ and $\varepsilon'<\varepsilon$,
the right-hand side above is no larger than
$\varepsilon\bar c(1+20\|f\|_\infty)$, which establishes (\ref{p3}).

To prove \hyperlink{equationP2}{(P2)}, it now suffices to
prove that for all $0<\delta<A<\infty$ there exists
$\zeta>0$ such that
%
%
\begin{equation}\label{p2}
\limsup_{M\to\infty} \sup_{\theta\in B(A/M)\setminus
B(\delta/M)} \phi_M(\theta) \le1-\zeta.
\end{equation}
Let $\tilde\phi_M(\theta) = \phi_M(\theta/M)$. By (iii),
$\tilde\phi_M(\theta)\to\tilde\phi(\theta)$ as $M\to\infty$,
and the convergence is uniform on compact sets. Since the
probability distribution with density $c_0f(x)$ on $B(1/2)$ is
not degenerate or of lattice type, $|\tilde\phi(\theta)|$ must be
bounded away from 1 on any compact set not containing 0. For
$0<\delta<A$, we may choose $\zeta>0$ such that
$\tilde\phi(\theta)< 1-\zeta$ for all $\theta\in
B(A)\setminus B(\delta)$. The uniform convergence
$\tilde\phi_M\to\tilde\phi$ on $B(A)\setminus B(\delta)$
now implies (\ref{p2}).

\section[Proof of Theorem 1.7]{Proof of Theorem \protect\ref{thm:uniform}}\label{sec4}
We continue to write $M$ for $M_L$. It
suffices to prove that
%
%
\begin{equation}\label{g1}
\lim_{L\to\infty}\sup_{x\in\mathbb{T}_L} L^2|
P_0(X^L_{s_L}=x) -L^{-2}| = 0.
\end{equation}
By pulling out the $y=0$ term from (\ref{ptL}), we see that
\[
L^2 | P_0(X^L_{s_L}=x) - L^{-2}|
= \biggl|\sum_{y\in\mathbb{T}'_L}
\phi^{s_L}_M(2\pi y/L) e^{2\pi i x y/L} \biggr|
\le\sum_{y\in\mathbb{T}'_L}
\phi^{s_L}_M(2\pi y/L).
\]
The limit (\ref{g1}) will follow from showing the
last sum tends to zero as
$L\to\infty$.

By \hyperlink{equationP1}{(P1)} there exists $\delta>0$
such that for large $L$,
%
%
\begin{equation}
1-\phi_M(2\pi y/L) \ge
\pi^2\sigma^2M^2|y|^2/L^2
\qquad\mbox{for all } y\in\mathbb{T}'_{\delta L/M}.
\end{equation}
This implies [recall (\ref{phit})] that
\[
\sum_{y\in\mathbb{T}'_{\delta L/M}}
\phi_M^{s_L}(2\pi y/L)
\le\sum_{y\in\mathbb{T}'_{\delta L/M}}
\exp(-s_L \pi^2\sigma^2M^2|y|^2/L^2).
\]
This last sum tends to 0 as $L\to\infty$ by comparison with
\[
\int_0^\infty e^{-\pi^2\sigma^2 s_L(M^2/L^2) r^2} r\, dr =
\frac{1}{2\pi^2\sigma^2s_LM^2/L^2} \to0
\]
since $s_LM^2/L^2\to\infty$ by assumption.

By \hyperlink{equationP2}{(P2)} and \hyperlink{equationP3}{(P3)}, there exists $\zeta>0$
such that for all large $L$,
\[
1-\phi_M(2\pi y /L) \ge\zeta\qquad\mbox{for all }
y\in\mathbb{T}_L\setminus\mathbb{T}_{\delta M/L}.
\]
This bound implies
\[
\sum_{y\in\mathbb{T}_{L}\setminus\mathbb{T}_{\delta L/M}}
\phi_M^{s_L}(2\pi y/L)
\le L^2\exp(- \zeta s_L )
\to0
\]
since $s_L/\log L \to\infty$ by assumption.
This completes the proof of (\ref{g1}).

\section[Proof of Theorem 1.2]{Proof of Theorem \protect\ref{thm:rick}}\label{sec5}
We continue to write $M$ for $M_L$.
To prove (\ref{eqn:exp-law-1}), it suffices in view of
(\ref{FL}) to establish the following facts:
%
%
\begin{equation}\label{eqn:GL0}
\lim_{L\to\infty} G_L(0,\lambda/L^2) =
\lambda^{-1} + 1
\end{equation}
and
\begin{equation}
\label{eqn:GLx}
{\lim_{L\to\infty} \sup_{x\in\mathbb{T}'_L}}
|G_L(x,\lambda/L^2) - \lambda^{-1} |
=0.
\end{equation}

\begin{pf*}{Proof of \protect(\ref{eqn:GL02})}
By (\ref{GL}),
%
%
\begin{equation}\label{goal0}
G_L(0,\lambda/L^2) = \lambda^{-1} +
\frac1{L^2 }\sum_{y\in\mathbb{T}'_L}
\frac{1}{1+\lambda/L^2 - \phi_M(2\pi
y/L)}
\end{equation}
and thus (\ref{eqn:GL0}) will follow from
%
%
\begin{equation}\label{goal1}
\lim_{L\to\infty}\frac1{L^2}\sum_{y\in\mathbb{T}'_L}
\frac{1}{1 - \phi_M(2\pi
y/L)} = 1.
\end{equation}

We will prove (\ref{goal1}) by breaking $\mathbb{T}'_L$ into
regions appropriate for utilizing \mbox{\hyperlink{equationP1}{(P1)}--\hyperlink{equationP3}{(P3)}}.
To prepare for this, fix $\varepsilon>0$. By
\hyperlink{equationP1}{(P1)}, there exists $\delta>0$
such that for all large $L$,
%
%
\begin{equation}\label{1-a}
\frac{1}{1-\phi_M(2\pi y/L)} \le\frac{1}
{\pi^2\sigma^2M^2|y|^2/L^2}
\qquad\mbox{for } y\in\mathbb{T}'_{\delta L/M}.
\end{equation}
By \hyperlink{equationP2}{(P2)} there exists $\delta'>0$ and $\zeta>0$ such that
for all large $L$,
%
%
\begin{equation}\label{2-a}
\frac{1}{1-\phi_M(2\pi y/L)} <
1/\zeta\qquad
\mbox{for } y\in\mathbb{T}_{\delta'L}\setminus\mathbb{T}_{\delta
L/M}.
\end{equation}
By \hyperlink{equationP3}{(P3)}, for any $0<a<\delta'$ and all large $L$,
%
%
\begin{equation}\label{3-a}
\biggl|\frac{1}{1-\phi_M(2\pi y/L)} -1 \biggr|<
\varepsilon\qquad
\mbox{for } y\in\mathbb{T}_{L}\setminus\mathbb{T}_{aL}.
\end{equation}

We claim that
%
%
\begin{eqnarray}\label{sum-1}
\lim_{L\to\infty} \frac1{L^2}\sum_{y\in
\mathbb{T}'_{\delta L/M}}
\frac{1}{1 - \phi_M(2\pi y/L)}
&=& 0,
\\
\label{sum-2}
\limsup_{L\to\infty}\frac{1}{L^2}\sum_{y\in\mathbb
{T}_{aL}\setminus\mathbb{T}
_{\delta L/M}}
\frac{1}{1 - \phi_M(2\pi y/L)}
&\le& a^2/\zeta
\end{eqnarray}
and
\begin{equation}
\label{sum-3}
\limsup_{L\to\infty} \biggl|
\frac1{L^2}\sum_{y\in\mathbb{T}_L\setminus\mathbb{T}_{aL}}
\frac{1}{1 - \phi_M(2\pi y/L)} - (1-a^2) \biggr|
\le\varepsilon.
\end{equation}
The bounds (\ref{sum-2}) and (\ref{sum-3}) are
immediate from (\ref{2-a}) and (\ref{3-a}).
For (\ref{sum-1}), we note that since
$M^2/\log L\to\infty$, (\ref{twopi}) implies
that
%
%
\begin{equation}\label{lim-1}
\lim_{L\to\infty}\frac{1}{M^2}
\sum_{y\in\mathbb{T}'_{L}}
\frac{1}{|y|^2} = 0.
\end{equation}
This fact and (\ref{1-a}) easily imply (\ref{sum-1}).
We note for later use that neither (\ref{sum-2}) nor
(\ref{sum-3}) require
$M^2/\log L\to\infty$, they hold for any $M\to\infty$ and
$\phi_M$ satisfying \hyperlink{equationP2}{(P2)} and \hyperlink{equationP3}{(P3)}.

Having established
(\ref{sum-1})--(\ref{sum-3}), we combine them to obtain
\[
\limsup_{L\to\infty} \biggl|
\frac1{L^2}\sum_{y\in\mathbb{T}'_L }
\frac{1}{1 - \phi_M(2\pi y/L)}-1 \biggr|
\le a^2/\zeta+ a^2
+\varepsilon.
\]
Let $a\downarrow0$ and then $\varepsilon\downarrow0$ to complete
the proof of (\ref{goal1}).
\end{pf*}
\begin{pf*}{Proof of \protect(\ref{eqn:GLx})}
After separating out the $y=0$ term as before, it suffices to
prove that
%
%
\begin{equation}\label{goal2}
\lim_{L\to\infty}\sup_{x\in\mathbb{T}'_L}
\frac1{L^2} \biggl|\sum_{y\in\mathbb{T}'_{L}}
\frac{e^{ 2\pi ix y/L}}{1-\phi_M(2\pi y/L)} \biggr|
= 0.
\end{equation}

In view of (\ref{sum-1}) and (\ref{sum-2}), we may concentrate
on the region $\mathbb{T}_L\setminus\mathbb{T}_{aL}$. By (\ref{3-a}),
uniformly in $x\in\mathbb{T}'_L$,
\[
\limsup_{L\to\infty} \frac1{L^2} \biggl|\sum_{y\in\mathbb
{T}_L\setminus
\mathbb{T}_{aL}}
e^{2\pi ixy/L} \biggl(
\frac{1}{1-\phi_M(2\pi y/L)} - 1 \biggr) \biggr|
\le\varepsilon.
\]
It is here we make use of (\ref{zero}). It implies that for
all $x\in\mathbb{T}'_L$,
\[
\biggl|\frac1{L^2}\sum_{y\in\mathbb{T}_L\setminus\mathbb
{T}_{aL}}e^{2\pi
ixy/L} \biggr|
= \biggl|-\frac1{L^2} \sum_{y\in\mathbb{T}_{aL}} e^{2\pi ixy/L} \biggr|
\le a^2.
\]
By the last two facts,
%
%
\begin{equation}\label{sum-4}
\limsup_{L\to\infty} \sup_{x\in\mathbb{T}'_L}\frac1{L^2} \biggl|\sum
_{y\in
\mathbb{T}_L\setminus
\mathbb{T}_{aL}}
\frac{e^{2\pi ixy/L}}{1-\phi_M(2\pi y/L)} \biggr|
\le\varepsilon+ a^2
\end{equation}
and we note here that (\ref{sum-4}) does not require that
$M^2/\log L\to\infty$.
Taken together, (\ref{sum-1}), (\ref{sum-2}) and (\ref{sum-4})
imply
%
%
\begin{equation}\label{sum-5}
\limsup_{L\to\infty} \sup_{x\in\mathbb{T}'_L}\frac1{L^2} \biggl|\sum
_{y\in
\mathbb{T}'_L}
\frac{e^{2\pi ixy/L}}{1-\phi_M(2\pi y/L)} \biggr|
\le\varepsilon+ a^2(1+1/\zeta).
\end{equation}
Let $a\to0$ and then $\varepsilon\to0$ to complete the proof of
(\ref{goal2}).
\end{pf*}
\begin{pf*}{Proof of \protect(\ref{eqn:mean1})}
By standard monotonicity
arguments,
%
%
\begin{equation}\label{unif-1}
P_x(H_L> u L^2) \to e^{-u} \qquad\mbox{uniformly in }
x\in\mathbb{T}'_L, u\ge0
\end{equation}
as $L\to\infty$. In particular, for all large $L$,
\[
P_x( H_L>L^2) \le e^{-1/2}\qquad
\mbox{for all }x\in\mathbb{T}'_L.
\]
By this bound and the Markov property,
\begin{eqnarray*}
P_x( H_L> kL^2) &=&
\sum_{y\in\mathbb{T}'_L}P_x\bigl( X^L_{(k-1)L}=y, H_L>(k-1)L^2\bigr)
P_y(H_L> L^2)\\
&\le& e^{-1/2} P_x\bigl( H_L> (k-1)L^2\bigr).
\end{eqnarray*}
Consequently, for all large $L$, $P_x(H_{L}> kL^2)\le
e^{-k/2}$ for $k\ge1$ and $x\in\mathbb{T}'_L$. This fact and
(\ref{unif-1}) easily imply (\ref{eqn:mean1}).
\end{pf*}

\section[Proof of Theorem 1.3]{Proof of Theorem \protect\ref{thm:main}}\label{sec6}
We continue to write $M$ for $M_L$.
The limit (\ref{eqn:exp-law}) follows easily from a little
algebra and the following
analogues of (\ref{eqn:GL0}), (\ref{eqn:GLx}):
%
%
\begin{equation}\label{eqn:GL02}
\lim_{L\to\infty} \frac{G_L(0,\lambda/L^2t_L)}{t_L} =
\lambda^{-1} + \rho+\frac1{\pi\sigma^2}
\end{equation}
and
\begin{equation}
\label{eqn:GLx2}
\lim_{L\to\infty} \sup_{x\in\mathcal{A}(\alpha,v_L)}
\biggl|\frac{G_L(x,\lambda/L^2t_L)}{t_L}
- \biggl[ \lambda^{-1}
+ \biggl(\frac{1-\alpha}{\pi\sigma^2} \biggr) \biggr] \biggr|
=0.
\end{equation}
The proofs of (\ref{eqn:GL02}) and (\ref{eqn:GLx2}) are
similar to the proofs of (\ref{eqn:GL0}) and
(\ref{eqn:GLx}), but require a bit more care.

Fix $\varepsilon>0$. By \hyperlink{equationP1}{(P1)}
there exist $\delta>0$
and functions $\psi_L$ such that
$\|\psi_L\|_\infty<\varepsilon$ and for all large $L$,
%
%
\begin{equation}\label{1-b}
\frac{1}{1-\phi_M(2\pi y/L)} = \frac{1+\psi_L(y)}
{2\pi^2\sigma^2M^2|y|^2/L^2}
\qquad\mbox{for } y\in\mathbb{T}'_{\delta L/M}.
\end{equation}
As before, we assume $\delta',\zeta>0$
are such that for all $0<a<\delta'$, (\ref{2-a})
and (\ref{3-a}) hold. Recall that
we are now assuming $M^2/\log L\to\rho<\infty$.
\begin{pf*}{Proof of \protect(\ref{eqn:GL02})}
The $y=0$ term in the sum for $G_L(0,\lambda/t_L)$
yields $\lambda^{-1}$, so it suffices to prove
that
%
%
\begin{equation}\label{goal3}
\lim_{L\to\infty}\frac{1}{L^2t_L} \sum_{y\in\mathbb
{T}'_L}\frac
{1}{1-\phi_M(2\pi
y/L)}
= \rho+\frac{1}{\pi\sigma^2}.
\end{equation}
We claim that:
%
%
\begin{eqnarray}\label{n1}
\limsup_{L\to\infty} \biggl|\frac{1}{L^2 t_L}\sum_{y\in\mathbb
{T}'_{\delta L/M}}
\frac{1}{1-\phi_M(2\pi y/L)} -
\frac{1}{\pi\sigma^2} \biggr|
&\le&\frac{\varepsilon}{\pi\sigma^2},
\\
\label{n2}
\limsup_{L\to\infty}\frac{1}{L^2t_L}\sum_{y\in\mathbb
{T}_{aL}\setminus
\mathbb{T}_{\delta L/M}}
\frac{1}{1 - \phi_M(2\pi y/L)}
&\le&\rho a^2/\zeta
\end{eqnarray}
and
\begin{equation}
\label{n3}
\limsup_{L\to\infty} \biggl|
\frac1{L^2t_L}\sum_{y\in\mathbb{T}_L\setminus\mathbb{T}_{aL}}
\frac{1}{1 - \phi_M(2\pi y/L)} - (1-a^2)\rho\biggr|
\le\varepsilon\rho.
\end{equation}
The limits (\ref{sum-2}) and (\ref{sum-3}) and the fact that
$1/t_L\to\rho$ imply (\ref{n2}) and (\ref{n3}), so consider
the region the region $\mathbb{T}'_{\delta L/M}$. By (\ref{1-b}),
\[
\frac{1}{L^2 t_L}\sum_{y\in\mathbb{T}'_{\delta L/M}}
\frac{1}{1-\phi_M(2\pi y/L)} =
\frac{1}{\log L}\sum_{y\in\mathbb{T}'_{\delta L/M}}
\frac{1+\psi_L(y)}{2\pi^2\sigma^2|y|^2}.
\]
By using (\ref{twopi}) above, we obtain (\ref{n1}).

Combining (\ref{n1})--(\ref{n3}) gives
\[
\limsup_{L\to\infty} \biggl|\frac{1}{L^2t_L}\sum_{y\in\mathbb{T}'_{L}}
\frac{1}{1 - \phi_M(2\pi y/L)} -\beta\biggr|
\le\frac{\varepsilon}{\pi\sigma^2} +
\rho a^2/\eta+ \rho(a^2+\varepsilon).
\]
Let $a\to0$ and then $\varepsilon\to0$ to complete the proof of
(\ref{goal3}).
\end{pf*}
\begin{pf*}{Proof of \protect(\ref{eqn:GLx2})}
Fix $0<\alpha<1$. (We will not give the slight changes in
proof needed to handle the cases $\alpha=0,1$.) It
suffices to prove that uniformly in $x\in\mathcal{A}(\alpha,v_L)$,
%
%
\begin{equation}\label{alphagoal}
\frac{1}{L^2t_L}
\sum_{y\in\mathbb{T}'_L} \frac{e^{2\pi ixy/L}}
{1-\phi_M(2\pi y/L)} \to\frac{1-\alpha}{\pi\sigma^2}\qquad
\mbox{as }L\to\infty.
\end{equation}
With $\varepsilon,\delta$ as before, we claim that
%
%
\begin{equation}\label{o1}\quad
\limsup_{L\to\infty}\sup_{x\in\mathcal{A}(\alpha,v_L)} \biggl|\frac{1}{L^2
t_L}\sum_{y\in\mathbb{T}'_{\delta L/M}}
\frac{e^{2\pi ixy/L}}{1-\phi_M(2\pi y/L)} -
\frac{1-\alpha}{\pi\sigma^2} \biggr|
\le\frac{\varepsilon}{\pi\sigma^2}.
\end{equation}
Given this, (\ref{sum-4}) and (\ref{n2}) imply
\begin{eqnarray*}
&&\limsup_{L\to\infty}\sup_{x\in\mathcal{A}(\alpha,v_L)} \biggl|\frac{1}{L^2
t_L}\sum_{y\in\mathbb{T}'_{L}}
\frac{e^{2\pi ixy/L}}{1-\phi_M(2\pi y/L)} -
\frac{1-\alpha}{\pi\sigma^2} \biggr|
\\
&&\qquad\le\frac{\varepsilon}{\pi\sigma^2}+\rho(\varepsilon
+a^2+a^2/\zeta),
\end{eqnarray*}
which is enough to establish (\ref{alphagoal}).

The first step in proving (\ref{o1}) is to use (\ref{1-b}) to obtain
%
%
\begin{equation}\label{f0}\hspace*{32pt}
\frac{1}{L^2 t_L}\sum_{y\in\mathbb{T}'_{\delta L/M}}
\frac{e^{2\pi ixy/L}}
{1-\phi_M(2\pi y/L)} =
\frac{1}{\log L}\sum_{y\in\mathbb{T}'_{\delta L/M}}
\frac{e^{2\pi ixy/L}}
{2\pi\sigma^2|y|^2}\bigl(1+\psi_L(y)\bigr).
\end{equation}
Next, we may replace $\mathbb{T}'_{\delta L/M}$ in
the right-hand side above with
$D'_{\delta L/M}$ because
%
%
\begin{equation}\label{f1}
\lim_{L\to\infty}\frac{1}{\log L}\sum_{y\in\mathbb{T}'_{\delta
L/M}\setminus
D'_{\delta L/M}} \frac{1}
{|y|^2} = 0
\end{equation}
by (\ref{2pilog2}).
Now, we break $\mathcal{A}_L(\alpha,v_L)$ into
the union of the smaller regions
\[
\mathcal{D}_L(\alpha,m) = D_{L^\alpha(\log L)^{m+1}}\setminus
D_{L^\alpha(\log L)^{m}},\qquad
m\in[-k,k)\cap\mathbb{Z}.
\]
We will prove that for each fixed $m$,
%
%
\begin{equation}\label{f2}
\lim_{L\to\infty}\sup_{x\in\mathcal{D}_L(\alpha,m)}
\biggl|
\frac{1}{\log L}\sum_{y\in D'_{\delta M/L}}
\frac{e^{2\pi i xy y/L}}{2\pi^2\sigma^2|y|^2} -
\frac{1-\alpha}{\pi\sigma^2} \biggr|
=0.
\end{equation}
Since (\ref{o1}) will follow
from (\ref{f0})--(\ref{f2}),
the problem now is to prove (\ref{f2}).

To do this, fix $m\in\mathbb{Z}$, let $K_L=L^{1-\alpha}(\log
L)^{-(m+1/2)}$, and consider the
regions $D_{\delta L/M}\setminus D_{K_L}$ and $D'_{K_L}$.
The bound (\ref{DKJbound}) implies that for all $x\in
\mathcal{D}_L(\alpha,m)$,
%
%
\begin{eqnarray}\label{f3}\qquad
\frac{1}{\log L}
\biggl|\sum_{y\in D_{\delta L/M}\setminus D_{K_L}}
\frac{e^{2\pi ixy/L}}{|y|^2}
\biggr|
&\le&
\frac{C_0}{
(\log L)(1\wedge K_L|2\pi x/L|) }
\nonumber\\
&\le&\frac{C_0}{\log L}\vee
\frac{C_0}{2\pi K_L(\log L)^{m+1}L^{\alpha-1}}
\to0 \\
\eqntext{\mbox{as }L\to\infty.}
\end{eqnarray}

To handle the sum over $D'_{K_L}$,
we make use of the fact that $e^{2\pi ixy/L}\approx1$
there. More precisely,
for $x\in\mathcal{D}_L(\alpha,m)$,
\begin{eqnarray*}
\frac{1}{\log L} \biggl|\sum_{y\in D'_{K_L}} \frac{e^{2\pi ixy/L}-1}{|y|^2}
\biggr|&\le&
\frac{1}{\log L}\sum_{y\in D'_{K_L}} \frac{2\pi|x|/L}{|y|}\\
&\le&2\pi L^{\alpha-1}(\log L)^{m} \sum_{y\in D'_{K_L}}\frac1{|y|}.
\end{eqnarray*}
Comparison with an integral shows there is a constant $C<\infty$ such
that
$\sum_{y\in D'_{K_L}}|y|^{-1} \le C K_L$, so it follows that
%
%
\begin{equation}\label{f4}
\lim_{L\to\infty} \sup_{x\in\mathcal{D}_L(\alpha,m)}
\frac{1}{\log L} \biggl|\sum_{y\in D'_{K_L}} \frac{e^{2\pi ixy/L}-1}{|y|^2}
\biggr| = 0.
\end{equation}
Coming to the main term at last, by (\ref{twopi}) we see that
%
%
\begin{eqnarray}\label{f5}\qquad
&&\frac{1}{2\pi^2\sigma^2\log L}\sum_{y\in D'_{K_L}} \frac{1}{|y|^2}
=\frac{\log K_L}{2\pi^2\sigma^2\log L}
\frac{1}{\log K_L}\sum_{y\in D'_{K_L}} \frac{1}{|y|^2}
\to\frac{1-\alpha}{\pi\sigma^2} \nonumber\\[-8pt]\\[-8pt]
&&\eqntext{\mbox{as }L\to\infty.}
\end{eqnarray}
Taken together, (\ref{f3})--(\ref{f5}) establish (\ref{f2}), as
required.
\end{pf*}
\begin{pf*}{Proof of \protect(\ref{eqn:mean})}
We proceed as in the proof of (\ref{eqn:mean1}) with just
a few changes. First, by (\ref{eqn:exp-law}) with
$\alpha=m=1$, there exists a finite $L_0$ such that for all
$L\ge L_0$, $P_y(H_L>L^2t_L)\le e^{-1/2\beta}$ for all
$y\in\mathbb{T}_L\setminus\mathbb{T}_{L/\log L}$. Next, by
Theorem \ref{thm:uniform}, there exists finite $L_1\ge
L_0$ such that for $L\ge L_1$ and all $x,y\in\mathbb{T}_L$,
$P_x(X^L_{L^2 t_L}=y) \le2/L^2 $. Therefore, for all $L\ge
L_1$ and $x\in\mathbb{T}'_L$,
\begin{eqnarray*}
P_x(H_L>2L^2t_L ) &\le&
P_x(X^L_{L^2t_L\in\mathbb{T}_{L/\log L}}) +
\sup_{y\in\mathbb{T}_L\setminus T_{L/\log L}}P_y(H_L>L^2t_L )\\
&\le&2|\mathbb{T}_{L/\log L}|/L^2 + e^{-1/2\beta}
\le2/(\log L)^2 + e^{-1/2\beta}.
\end{eqnarray*}
It follows that for some finite $L_2\ge L_1$, if $L\ge L_2$
then
\[
\sup_{x\in\mathbb{T}'_L} P_x( H_L>2L^2t_L) \le e^{-1/3\beta}.
\]
Iterating as in the proof of (\ref{eqn:mean1}), we obtain
%
%
\begin{equation}\label{unifbnd}
\sup_{x\in\mathbb{T}'_L} P_x( H_L>2kL^2t_L) \le e^{-k/3\beta}
\end{equation}
for all $L\ge L_2$.

Now for a fixed $0\le\alpha\le1$ and $k>0$,
(\ref{eqn:exp-law}) implies
%
%
\begin{eqnarray}
P_x(H_L>uL^2t_L) \to(1-q)e^{-u/\beta}
\nonumber\\[-8pt]\\[-8pt]
\eqntext{\mbox{uniformly in } x\in\mathcal{A}(\alpha,v_L), u\ge0,}
\end{eqnarray}
as $L\to\infty$. The limit (\ref{eqn:mean}) is a consequence
of this fact and (\ref{unifbnd}).
\end{pf*}

\section[Example 1.5]{Example \protect\ref{myexample}}\label{sec7}
In this section, we verify the claims made in
Example \ref{myexample}. We first check that
\begin{eqnarray*}
\frac{\sigma^2_{M_L}}{M^2_L} &=& \frac{1}{M_L^2}\sum_x x_1^2
q_{M_L}(x) \\
&=& \frac{c}{M_L^2}\sum_x x_1^2 u_{M_L}(x) +
\frac{1-c}{M_L^2}\sum_x x_1^2 q_{0}(x)
\to c\int_{B(1/2)} x^2_1 \,dx \\
&=& \frac c{12}\qquad
\mbox{as }L\to\infty,
\end{eqnarray*}
so (\ref{insuffconds}) holds with
$\sigma^2=c/12$. We turn now to the proof of
(\ref{eqn:mix}).

Let $\hat u_{M_L}(\theta) = \sum_x u_{M_L}(x)e^{i\theta x}$. Our
first step is to establish the analogues of
\hyperlink{equationP1}{(P1)}--\hyperlink{equationP3}{(P3)} for $\phi_{M_L}(\theta) =
c\hat u_{M_L}(\theta) + (1-c)\hat q_0(\theta)$. By Proposition
\ref{prop:suff},\vspace*{1pt} $\hat u_{M_L}$ satisfies
\hyperlink{equationP1}{(P1)}--\hyperlink{equationP3}{(P3)} with $\sigma^2= 1/12$.
Furthermore, it is easy to check that $\hat q_{M_0}$ satisfies: for
all $\varepsilon>0$ there exists $\delta>0$ such that
\[
\frac{1-\hat q_{0}(\theta)}{\sigma_{0}^2|\theta|^2/2}\in
(1-\varepsilon,1+\varepsilon)\qquad
\mbox{for all }\theta\in B'(\delta).
\]
With this it is easy to see that
the following versions of \hyperlink{equationP1}{(P1)}--\hyperlink{equationP3}{(P3)} hold for $\phi_{M_L}$.
\begin{longlist}[(P1)$'$]
\item[(P1)$'$]\hypertarget{equationP1prime}
For $\varepsilon>0$ there exists
$\delta>0$ such that for all large $L$,
\[
\frac{1}{1-\phi_{M_L}(2\pi y/L)} =
\frac{1+\psi_L(y)}{cM_L^2\pi^2|y|^2/6L^2}
\qquad\mbox{for all }y\in\mathbb{T}'_{\delta L/M_L},
\]
where $\|\psi_L\|_\infty\le\varepsilon$.
\item[(P2)$'$]\hypertarget{equationP2prime}
For $\delta>0$ there exists $\delta'>0$ and
$\zeta>0$ such that for all large $L$,
\[
1-\phi_{M_L}(2\pi y/L) \ge c \zeta\qquad
\mbox{for all }y\in\mathbb{T}_{\delta' L}\setminus\mathbb
{T}_{\delta L/M_L}
.
\]
\item[(P3)$'$]\hypertarget{equationP3prime}
For fixed $0<a<1$,
\[
\lim_{L\to\infty}\sup_{y\in\mathbb{T}_L\setminus\mathbb
{T}_{aL}} \biggl|
\frac{1}{1-\phi_{M_L}(2\pi y/L)} -
\frac{1}{c+(1-c)(1-\hat q_0(2\pi y/L))}
\biggr| = 0.
\]
\end{longlist}

With the above in place, the next step is to prove that
\[
\lim_{L\to\infty} G_L(0,\lambda/L^2) = \lambda^{-1}+\beta_0
\]
or equivalently
%
%
\begin{equation}\label{GL0-mix}
\lim_{L\to\infty}\frac1{L^2}\sum_{
y\in\mathbb{T}_L'} \frac{1}{1-\phi_{M_L}(2\pi y/L)}= \beta_0.
\end{equation}
To do this fix $\varepsilon>0$, choose $\delta,\delta'$ as in
\hyperlink{equationP1prime}{(P1)$'$} and \hyperlink{equationP2prime}{(P2)$'$},
and break $\mathbb{T}'_L$ into the usual subregions.

Applying \hyperlink{equationP1prime}{(P1)$'$}, we have
\[
\frac1{L^2}\sum_{y \in\mathbb{T}'_{\delta L/M_L}}
\frac{1}{1 - \phi_{M_L}(2\pi y/L)}
=
\frac{6}{cM^2_L\pi^2}\sum_{y \in\mathbb{T}'_{\delta L/M_L}}
\frac{1+\psi_L(y)}{|y|^2}.
\]
This implies, using (\ref{twopi}),
%
%
\begin{equation}\label{GL0-mix1}
\limsup_{L\to\infty} \biggl|\frac1{L^2}\sum_{y \in\mathbb{T}'_{\delta L/M_L}}
\frac{1}{1 - \phi_{M_L}(2\pi y/L)} -
\frac{12}{c\pi} \biggr| \le\frac{12\varepsilon}{c\pi},
\end{equation}
where we have used $M^2_L/\log L\to1$.
Next, for $0<a<\delta'$, \hyperlink{equationP2prime}{(P2)$'$} implies
%
%
\begin{equation}\label{GL0-mix2}
\limsup_{L\to\infty}\frac{1}{L^2}\sum_{y\in\mathbb
{T}_{aL}\setminus
\mathbb{T}_{\delta L/M_L}}
\frac{1}{1 - \phi_{M_L}(2\pi y/L)}
\le\frac{a^2}{c\zeta}.
\end{equation}
By \hyperlink{equationP3prime}{(P3)$'$} and continuity,
\begin{eqnarray*}
\frac1{L^2}\sum_{y\in\mathbb{T}_L\setminus\mathbb{T}_{aL}}
\frac{1}{1-\phi_{M_L}(2\pi y/L)} &\to&
\int_{B(1/2)\setminus B(a/2)}
\frac{d\theta}{c+ (1-c) (1-\hat q_0(2\pi\theta))}\\
&=&\frac{1}{(2\pi)^{2}}\int_{B(\pi)\setminus B(a\pi)}
\frac{d\theta}{1-(1-c)\hat q_0(\theta)}.
\end{eqnarray*}
Let $a\downarrow0$ and then $\varepsilon\downarrow0$ in
(\ref{GL0-mix1}) and (\ref{GL0-mix2})
to complete the proof of (\ref{GL0-mix}).

The final task is to prove that
\[
{\lim_{L\to\infty}\sup_{x\in\mathbb{T}_L\setminus\mathbb
{T}_{\ell_L}}}
|G_L(x,\lambda/L^2) - \lambda^{-1}| = 0
\]
or equivalently
%
%
\begin{equation}\label{GLx-mix}
\lim_{L\to\infty}\sup_{x\in\mathbb{T}_L\setminus\mathbb
{T}_{\ell_L}} \biggl|
\frac1{L^2}\sum_{y\in\mathbb{T}'_L}
\frac{e^{2\pi i x y/L}}{1-\phi_{M_L}(2\pi
y/L)} \biggr| = 0.
\end{equation}

Consider the region $\mathbb{T}'_{\delta L/M}$. By (\ref{2pilog2}), we
may replace $\mathbb{T}'_{\delta L/M}$ with $D'_{\delta L/M}$, at the
cost of a negligible error. We break $D'_{\delta L/M}$ into two pieces.
By (\ref{twopi}),
%
%
\begin{equation}\label{GLx-mix1}
\lim_{L\to\infty}\frac{1}{\log L}\sum_{y\in D'_{L^\varepsilon
}}\frac1{|y|^2}
= 2\pi\varepsilon.
\end{equation}
By (\ref{phiest}), for all $x\in\mathbb{T}_L\setminus\mathbb
{T}_{\ell}$
%
%
\begin{eqnarray}\label{GLx-mix2}\hspace*{32pt}
\biggl|\frac{1}{\log L}\sum_{y\in D_{\delta L/M}\setminus D_{L^\varepsilon}}
\frac{e^{2\pi i x y/L}}{|y|^2} \biggr|
&\le&\frac{C_0}{\log L} \vee
\frac{C_0L}{K_L 2\pi|x|} \nonumber\\[-8pt]\\[-8pt]
&\le&
\frac{C_0}{\log L} \vee\frac{C_0 L^{1-\varepsilon}}{2\pi\ell_L}
\to0 \qquad\mbox{as } L\to\infty.\nonumber
\end{eqnarray}
By \hyperlink{equationP1prime}{(P1)$'$} and the above,
%
%
\begin{equation}
\limsup_{L\to\infty}\sup_{x\in\mathbb{T}_L\setminus\mathbb
{T}_{\ell_L}} \biggl|
\frac1{L^2}\sum_{y\in\mathbb{T}'_{\delta M_L/L}}
\frac{e^{2\pi i x y/L}}{1-\phi_{M_L}(2\pi
y/L)} \biggr| \le\frac{12\varepsilon}{c\pi}
\end{equation}
and combining this with (\ref{GL0-mix2}) gives
%
%
\begin{equation}\label{GLx-mix3}
\limsup_{L\to\infty}\sup_{x\in\mathbb{T}_L\setminus\mathbb
{T}_{\ell_L}} \biggl|
\frac1{L^2}\sum_{y\in\mathbb{T}'_{aL}}
\frac{e^{2\pi i x y/L}}{1-\phi_{M_L}(2\pi
y/L)} \biggr| \le\frac{12\varepsilon}{c\pi}
+ \frac{a^2}{c\zeta}.
\end{equation}

Now consider the region $\mathbb{T}_L\setminus\mathbb{T}_{aL}$.
By \hyperlink{equationP3prime}{(P3)$'$}, for all large $L$ and $x\in\mathbb{T}_L$,
%
%
\begin{equation}\label{GLx-mix4}
\frac{1}{L^2}\sum_{{y\in\mathbb{T}_L\setminus\mathbb{T}_{aL}}}
\biggl| \frac{ e^{2\pi i x y/L}}
{1-\phi_{M_L}(2\pi y/L)} - \frac{e^{2\pi i x
y/L}}{1-(1-c)\hat q_0(2\pi y/L)}
\biggr| \le\varepsilon.
\end{equation}
For integers $K>0$ define $\Gamma_{L,K} = \{ z\in K\mathbb{Z}^2\dvtx
z+\mathbb{T}_K\subset\mathbb{T}_L\setminus\mathbb{T}_{aL}\} $, and
note that
$|\Gamma_{L,K}|\le L^2/K^2$ and $|(\mathbb{T}_L\setminus
\mathbb{T}_{aL})\setminus\bigcup_{z\in\Gamma_{L,K}}(z+\mathbb{T}_K)
|\le
8LK$. By the trivial bound $1-(1-c)\hat q_0(\theta)\ge c$ and
(\ref{GLx-mix4}),
%
%
\begin{eqnarray}\label{GLx-mix5}
&&\biggl|\frac{1}{L^2}\sum_{{y\in\mathbb{T}_L\setminus\mathbb{T}_{aL}}}
\frac{ e^{2\pi i x y/L}}
{1-\phi_{M_L}(2\pi y/L)} \nonumber\\[-8pt]\\[-8pt]
&&\qquad{} - \frac{1}{L^2}\sum_{z\in\Gamma_{L,K}}\sum_{y\in z+\mathbb{T}_K}
\frac
{e^{2\pi i x
y/L}}{1-(1-c)\hat q_0(2\pi y/L)}
\biggr| \le\varepsilon+\frac{8K}{cL}.\nonumber
\end{eqnarray}

By the continuity of $\hat q_0$, there exists $\delta''>0$
such that if $\theta,\theta'\in
B(\pi)$ and $|\theta-\theta'|<\delta''$ then
\[
\biggl|
\frac1{1-(1-c)\hat q_0(\theta)} -
\frac1{1-(1-c)\hat q_0(\theta')} \biggr|
<\varepsilon.
\]
Assuming $K<\delta'' L$, this implies
%
%
\begin{eqnarray}\label{GLx-mix6}
&&\biggl|
\frac{1}{L^2} \sum_{z\in\Gamma_{L,K}}
\sum_{y\in{z+\mathbb{T}_K}}\frac{e^{2\pi i x y/L}}{1-(1-c)\hat
q_0(2\pi y/L)}
\nonumber\\[-8pt]\\[-8pt]
&&\qquad{} - \frac{1}{L^2} \sum_{z\in\Gamma_{L,K}}
\frac{e^{2\pi ixz/L}}{1-(1-c)\hat q_0(2\pi z/L)}
\sum_{y\in{\mathbb{T}_K}}e^{2\pi i x y/L}
\biggr| < \varepsilon.\nonumber
\end{eqnarray}
Now (\ref{phiest}) can be applied, giving
%
%
\begin{eqnarray}\label{GLx-mix7}
&&\biggl| \frac{1}{L^2} \sum_{z\in\Gamma_{L,K}}
\frac{e^{2\pi ixz/L}}{1-(1-c)\hat q_0(2\pi
z/L)}\sum_{y\in{\mathbb{T}_K}}e^{2\pi i x y/L}
\biggr| \nonumber\\[-8pt]\\[-8pt]
&&\qquad\le
\frac{|\Gamma_{L,K}|}{cL^2}
\biggl| \sum_{y\in{\mathbb{T}_K}}e^{2\pi i x y/L}
\biggr|
\le\frac{(4(K+1)(1+ L/2\pi\ell_L))}{cK^2}\nonumber
\end{eqnarray}
for all $x\in\mathbb{T}_L\setminus\mathbb{T}_{\ell_L}$.
Taken together (\ref{GLx-mix3}) and (\ref{GLx-mix5})--(\ref{GLx-mix7})
yield
\begin{eqnarray*}
&&\limsup_{L\to\infty}\sup_{x\in\mathbb{T}_L\setminus
\mathbb{T}_{\ell_L}} \biggl|\frac{1}{L^2}\sum_{{y\in\mathbb{T}'_L}}
\frac{ e^{2\pi i x y/L}}
{1-\phi_{M_L}(2\pi y/L)} \biggr| \\
&&\qquad\le
\frac{12\varepsilon}{c\pi}+
\frac{a^2}{c\zeta} + 2\varepsilon
+ \limsup_{L\to\infty} \biggl(\frac{8K}{cL} +
\frac{(4(K+1)(1+ L/2\pi\ell_L))}{K^2}\biggr).
\end{eqnarray*}
If set $K=L/\sqrt{\ell}$, then the limsup above is 0.
Let $a\downarrow0$ and the $\varepsilon\downarrow0$ to finish the
proof.

\section*{Acknowledgment}
It is a pleasure to thank Rick Durrett for suggesting this
problem.


%
\printaddresses

\end{document}